\documentclass[12pt]{article}

\usepackage{amsmath,amssymb,amsthm,amscd}
\usepackage[colorlinks, linkcolor=blue, anchorcolor=gree, citecolor=red]{hyperref}
\usepackage[numbers,sort&compress]{natbib}
\usepackage{graphicx}

\begin{document}

\newcommand{\Cyc}{{\rm{Cyc}}}

\newtheorem{thm}{Theorem}[section]
\newtheorem{pro}[thm]{Proposition}
\newtheorem{lem}[thm]{Lemma}
\newtheorem{cor}[thm]{Corollary}
\theoremstyle{definition}
\newtheorem{ex}[thm]{Example}
\newtheorem{remark}[thm]{Remark}
\newcommand{\bth}{\begin{thm}}
\renewcommand{\eth}{\end{thm}}
\newcommand{\bex}{\begin{examp}}
\newcommand{\eex}{\end{examp}}
\newcommand{\bre}{\begin{remark}}
\newcommand{\ere}{\end{remark}}

\newcommand{\bal}{\begin{aligned}}
\newcommand{\eal}{\end{aligned}}
\newcommand{\beq}{\begin{equation}}
\newcommand{\eeq}{\end{equation}}
\newcommand{\ben}{\begin{equation*}}
\newcommand{\een}{\end{equation*}}

\newcommand{\bpf}{\begin{proof}}
\newcommand{\epf}{\end{proof}}
\renewcommand{\thefootnote}{}

\def\beql#1{\begin{equation}\label{#1}}
\title{\Large\bf  Non-cyclic graphs of (non)orientable genus one}

\author{{\sc Xuanlong Ma
}\\[15pt]
{\small Sch. Math. Sci. {\rm \&} Lab. Math. Com. Sys.,}\\
{\small Beijing Normal University, 100875, Beijing, China.}\\
}

 \date{}

\maketitle

\begin{abstract}
Let $G$ be a finite non-cyclic group. The non-cyclic graph $\Gamma_G$ of $G$
is the graph whose vertex set is $G\setminus \Cyc(G)$, two distinct vertices being adjacent if they do not generate a cyclic subgroup, where $\Cyc(G)=\{a\in G: \langle a,b\rangle\ \text{is cyclic for each }  b\in G\}$.
In this paper, we classify all finite non-cyclic groups $G$ such that $\Gamma_G$ has  (non)orientable genus one.
\end{abstract}



\footnote{E-mail addresses: xuanlma@mail.bnu.edu.cn.}
\section{Introduction}
All graphs in this paper are undirected, with no loops or multiple
edges.
A graph $\Gamma$ is called a {\em planar graph} if $\Gamma$ can be drawn in the plane so that no two of its edges cross each other, and in this case we say that $\Gamma$ can be embedded in the plane.
For a non-planar graph, it can be  embedded in some surface obtained from the sphere by attaching some handles or crosscaps.
We denote by $\mathbb{S}_k$ a sphere with $k$ handles and by $\mathbb{N}_k$ a sphere with $k$ crosscaps.
Note that both $\mathbb{S}_0$ and $\mathbb{N}_0$ are the sphere itself, and $\mathbb{S}_1$
and $\mathbb{N}_1$ are the torus and the projective plan, respectively.
The smallest non-negative integer $k$
such that a graph $\Gamma$ can be embedded on $\mathbb{S}_k$ (resp. $\mathbb{N}_k$) is called the {\em orientable genus} or {\em genus} (resp. {\em nonorientable genus}) of $\Gamma$, and is denoted by $\gamma(\Gamma)$ (resp. $\overline{\gamma}(\Gamma)$).

The problem of finding the graph genus is NP-hard \cite{T89}. The (non)orientable genera of some graphs constructed from some algebraic structures have been studied, for instance, see
\cite{W06,CSW10,RD15,A15,D14}.

All groups considered in this paper are finite. Denote by $\mathbb{Z}_n$ and  $D_{2n}$ the cyclic group of order $n$ and the dihedral group of order $2n$, respectively.
Let $G$ be a non-cyclic group. The {\em cyclicizer} $\Cyc(G)$ of $G$ is
$$\{a\in G: \langle a,b\rangle\ \text{is cyclic for each }  b\in G\}.$$
and is a normal subgroup of $G$ (see \cite{OP92}).
The non-cyclic graph $\Gamma_G$ of $G$
is the graph whose vertex set is $G\setminus \Cyc(G)$, and two distinct vertices being adjacent if they do not generate a cyclic subgroup.
The non-cyclic graph $\Gamma_G$ was first considered by
Abdollahi and  Hassanabadi \cite{AH07} and they studied the properties of the graph and established some graph theoretical properties (such as regularity) of this graph in terms
of the group ones. In \cite{AH09}, Abdollahi and  Hassanabadi classified all non-cyclic groups $G$ such that $\Gamma_G$ is planar.

A natural question is the following: Which finite non-cyclic groups have their non-cyclic graphs have(non)orientable genus one ?
The goal of the paper is to find all non-cyclic graphs of (non)orientable genus one.
Our main results are the following theorems.

\begin{thm}\label{th1}
Let $G$ be a finite non-cyclic group. Then $\Gamma_{G}$ has genus one
if and only if $G$ is isomorphic to one of the following groups:
\begin{equation}\label{g1}
\mathbb{Z}_3^2,~\mathbb{Z}_2^3,~\mathbb{Z}_2\times \mathbb{Z}_4,~D_8,~\mathbb{Z}_2 \times \mathbb{Z}_6.
\end{equation}
\end{thm}

\begin{thm}\label{th2}
Let $G$ be a finite non-cyclic group. Then $\Gamma_{G}$ has nonorientable genus one
if and only if $G$ is isomorphic to $\mathbb{Z}_2\times \mathbb{Z}_4$ or $D_8$.
\end{thm}

\section{Preliminaries}\label{sec:}

An element of order $2$ in a
group is called an {\em involution}.
Let $G$ be a group and $g$ be an element of $G$.
Denote by $|G|$ and $|g|$ the orders of $G$ and  $g$, respectively.
We denote the symmetric group on $n$ letters and the quaternion group of order $8$ by $S_n$ and
$Q_8$, respectively.
Also $\mathbb{Z}_n^m$ is used for the $m$-fold direct product of the cyclic group $\mathbb{Z}_n$ with itself. In the following, we state some results which we need in the sequel.


\begin{lem}{\rm (\cite[Proposition 4.3]{AH09})}\label{planar}
$\Gamma_G$ is planar if and only if $G$ is isomorphic to $\mathbb{Z}_2\times \mathbb{Z}_2$, $S_3$ or $Q_8$.
\end{lem}

Let $\Gamma$ be a graph.
Denote by $V(\Gamma)$ and $E(\Gamma)$ the vertex set and the edge set of $\Gamma$, respectively.
We use the natation $\lceil x\rceil$ to denote the least integer that is greater than or
equal to $x$.
Denote by $K_n$ and $K_{m,n}$ the complete graph of order $n$ and the complete bipartite graph,
respectively.
The following result from \cite{Whi} gives the (non)orientable genus of a complete graph and a complete multipartite graph.

\begin{lem}{\rm (\cite{Whi})}\label{ccgenus}
Let $n$ be an integer at least $3$. Then

$(a)$ $\gamma(K_n)=\lceil \frac{1}{12}(n-3)(n-4)\rceil$.

$(b)$  $\overline{\gamma}(K_n)=\lceil \frac{1}{6}(n-3)(n-4)\rceil$ if $n\ne 7$ and $\overline{\gamma}(K_7)=3$.

$(c)$ $\gamma(K_{m,n})=\lceil \frac{1}{4}(m-2)(n-2)\rceil$.

$(d)$ $\overline{\gamma}(K_{m,n})=\lceil \frac{1}{2}(m-2)(n-2)\rceil$.

$(e)$ $\gamma(K_{n,n,n})=\frac{1}{2}(n-1)(n-2)$.

$(f)$ $\gamma(K_{n,n,n,n})=(n-1)^2$ for $n\ne 3$ and $\gamma(K_{3,3,3,3})=5$.
\end{lem}

\begin{lem}{\rm (\cite[pp. 252, Theorem 9.7.3]{S64})}\label{unpg}
Suppose that $G$ is a $p$-group for some prime $p$ and has a unique
subgroup of order $p$.
If $p=2$,  then $G$ is cyclic or generalized quaternion. If $p>2$, then $G$ is cyclic.
\end{lem}

The following result is one of Sylow theorems.

\begin{thm}\label{Sythm}
Suppose that $G$ is a group and $p$ is a prime divisor of $|G|$. Then the number of
Sylow $p$-subgroups is congruent to $1$ modulo $p$. In particular, the number of
subgroups of order $p$ is congruent to $1$ modulo $p$.
\end{thm}


Denote by $\varphi$ the Euler's totient function.

\begin{lem}\label{fact5}
Let $G$ be a non-cyclic group, $p,q$ two distinct primes and $m$ a positive integer at least $1$.
If $\gamma(\Gamma_G)=1$, the each of the following statements does not hold:

$(a)$ $G$ has $4$ cyclic subgroups of order $p^m$ and an element of order $q$, where
$\varphi(p^m)\ge 2$.

$(b)$ $G$ has $4$ cyclic subgroups of order $3$ and $|G|\ge 10$.

$(c)$ $G$ has $3$ cyclic subgroups of order $4$ and an element of order $q^m$, where
$\varphi(q^m)\ge 3$ and $q\ne 2$.

$(d)$ $G$ has $7$ cyclic subgroups of order $2$ and an element of order $q$, where
$q\ne 2$.

$(e)$ $G$ has $3$ cyclic subgroups of order $p^m$, where
$\varphi(p^m)\ge 4$.

$(f)$ $G$ has $2$ cyclic subgroups of order $p^m$, where
$\varphi(p^m)\ge 5$.

\end{lem}
\bpf $(a)$ Suppose, for a contradiction,
that $(a)$ holds. Let $\langle a\rangle$, $\langle b\rangle$, $\langle c\rangle$, $\langle d\rangle$ be
$4$ cyclic subgroups of order $p^m$ of $G$ and $g$ be an element of order $q$.
If $g$ and each element of $\{a,a^{-1},b,b^{-1},c,c^{-1},d,d^{-1}\}$ cannot generate a cyclic subgroup, then the induced subgraph by  $\{a,a^{-1},b,b^{-1},c,c^{-1},d,d^{-1},g\}$ has a subgraph isomorphic to $K_{4,5}$ that has partition sets
$\{c,c^{-1},d,d^{-1}\}$ and $\{a,a^{-1},b,b^{-1},g\}$ and so $\gamma(\Gamma_G)\ge \gamma(K_{4,5})=2$,  a contradiction. Thus, we may
suppose that $g$ and an element of $\{a,a^{-1},b,b^{-1},c,c^{-1},d,d^{-1}\}$ can generate a cyclic subgroup $\langle h\rangle$. Without loss of generality, let $\langle a,g\rangle=\langle h\rangle$. Then $h\in V(\Gamma_G)$ and thereby, the induced subgraph by  $\{a,a^{-1},b,b^{-1},c,c^{-1},d,d^{-1},h\}$ has a subgraph isomorphic to $K_{4,5}$ that has partition sets
$\{c,c^{-1},d,d^{-1}\}$ and $\{a,a^{-1},b,b^{-1},h\}$. So
$\gamma(\Gamma_G)\ge \gamma(K_{4,5})> 1$, also a contradiction.

It is similar to the proof of $(a)$, we can prove $(b)$, $(c)$ and $(d)$.

$(e)$ Assume, to the contrary, that $(e)$ holds.
Take $4$ generators in every cyclic subgroup of order $p^m$. Then it is easy to see that
the induced subgraph by the generators  has a subgraph isomorphic to $K_{4,4,4}$ that has
genus $3$ by Lemma \ref{ccgenus}, a contradiction.

$(f)$ It is similar to the proof of $(e)$.
\epf

\begin{lem}\label{p-group}
Let $G$ be a non-cyclic $p$-group, where $p$ is a prime. Then $\gamma(\Gamma_G)=1$ if and only if $G$ is isomorphic to one of the following groups:
\begin{equation}\label{g2}
\mathbb{Z}_3^2,~\mathbb{Z}_2^3,~\mathbb{Z}_2\times \mathbb{Z}_4,~D_8.
\end{equation}
\end{lem}
\bpf
Note that $\Gamma_{\mathbb{Z}_3^2}\cong K_{2,2,2,2}$, $\Gamma_{\mathbb{Z}_2^3}\cong K_{7}$ and each of $\Gamma_{\mathbb{Z}_2\times \mathbb{Z}_4}$ and $\Gamma_{D_8}$ is a subgraph of $K_7$.
By Lemma \ref{ccgenus}, we see that $\Gamma_G$ has genus one for each group $G$ in (\ref{g2}).
We next assume that $\gamma(\Gamma_G)=1$.

Suppose that $p\ge 3$.
Then, by Lemma \ref{unpg} and Theorem \ref{Sythm} we have that $G$ has at least $4$ subgroups of order $p$. It follows from $(e)$ and $(b)$ of Lemma \ref{fact5} that
$|G|\le 9$. This implies that $G\cong \mathbb{Z}_3^2$, as desired.

Now suppose that $p= 2$.
If $|G|\le 8$, then $G$ is isomorphic to $\mathbb{Z}_2^2$,
$Q_8$, $\mathbb{Z}_2^3$, $\mathbb{Z}_2\times \mathbb{Z}_4$ or $D_8$ and by Lemma \ref{planar}
we get the desired result.
Thus, we may suppose that $|G|\ge 16$. If $G$ is generalized quaternion, then $G$ has a subgroup
$\langle x\rangle$ of order $8$ and contains at least $4$ elements $y_1,\cdots,y_4$ of order $4$
that do not belong to $\langle x\rangle$, and so
$\Gamma_G$ has a subgraph isomorphic to $K_{6,4}$ that has partition sets
$\{x^i: 0<i<8, i\ne 4\}$ and $\{y_j: j=1,\ldots,4\}$, a contradiction by Lemma \ref{ccgenus}.
Therefore, by Lemma \ref{unpg} we may assume that  $G$ has at least $3$ involutions.

\medskip
\noindent {\bf Case 1.} $G$ has an element $g$ of order $8$.
\medskip

Suppose that $G$ has two distinct subgroups $\langle g\rangle$, $\langle h\rangle$ of order $8$. Then we may pick an involution $a$ in $G\setminus (\langle g\rangle\cup \langle h\rangle)$.
Now we get a subgraph of $\Gamma_G$ isomorphic to $K_{5,4}$ that has partition sets
$\{g,g^3,g^5,g^7,a\}$ and $\{h,h^3,h^5,h^7\}$, a contradiction by Lemma \ref{ccgenus}.

Thus, we may suppose that $G$ has a unique subgroup $\langle u\rangle$ of order $8$, which is normal in $G$. Take an involution $b$ that does not belong to $\langle u\rangle$. Then
$\langle u,b\rangle$ is a subgroup of order $16$ and has precisely one subgroup of order $8$.
Since $b\notin \langle u\rangle$, $\langle u,b\rangle$ is not cyclic. Note that $G$ is not  generalized quaternion.
By verifying the groups of order $16$, we get that $G\cong D_{16}$ or $QD_{16}$, where
$QD_{16}=\langle a,b: a^8=b^2=1,bab=a^3\rangle$.
If $G\cong D_{16}$, then $G$ has $9$ involutions which induce a subgraph
isomorphic to $K_{9}$, $\gamma(\Gamma_G)\ge \gamma(K_9)=3$ by Lemma \ref{ccgenus}, a contradiction. Note that $QD_{16}$ has only $6$ elements of order $4$.
If $G\cong QD_{16}$, then the subgraph induced by
$6$ elements of order $4$ and $4$ elements of order $8$ has
a subgraph isomorphic to $K_{6,4}$, also a contradiction.

\medskip
\noindent {\bf Case 2.} $G$ has no elements of order $8$.
\medskip

If $G$ has no elements of order $4$, then $\Gamma_G$ is isomorphic to $K_{|G|-1}$, a contradiction as $\gamma(K_{|G|-1})>1$ for $|G|\ge 16$.
Thus, in this case we may assume that $\pi_e(G)=\{1,2,4\}$. Note that $|G|\ge 16$.
Since all involutions induce a complete graph and $\gamma(K_8)\ge 2$, $G$ has at least
$8$ elements of order $4$.
Since a power graph induced by $10$ elements of order $4$ of a group has
a subgraph isomorphic to $K_{6,4}$ that has genus two,
$G$ has precisely $8$ elements of order $4$ and $7$ involutions.
Take an involution $a$ in $G$ that does not belong to any subgroup of order $4$. Then
it is easy to see that the subgraph induced by all elements of order $4$ and $a$ has
a subgraph isomorphic to $K_{5,4}$ that has genus two, a contradiction.
\epf

\section{Proof of the main theorems}

{\noindent \em Proof of Theorem \ref{th1}.}
Note that $\Gamma_{\mathbb{Z}_2 \times \mathbb{Z}_6}$ has a subgraph isomorphic to $K_{3,3}$. Hence $\gamma(\Gamma_{\mathbb{Z}_2 \times \mathbb{Z}_6})\ge 1$. On the other hand,
we can embed $\Gamma_{\mathbb{Z}_2 \times \mathbb{Z}_6}$ into the tours as shown in Figure
\ref{fn3}. This implies that $\gamma(\Gamma_{\mathbb{Z}_2 \times \mathbb{Z}_6})= 1$.
\begin{figure}[hptb]
  \centering
  \includegraphics[width=8cm]{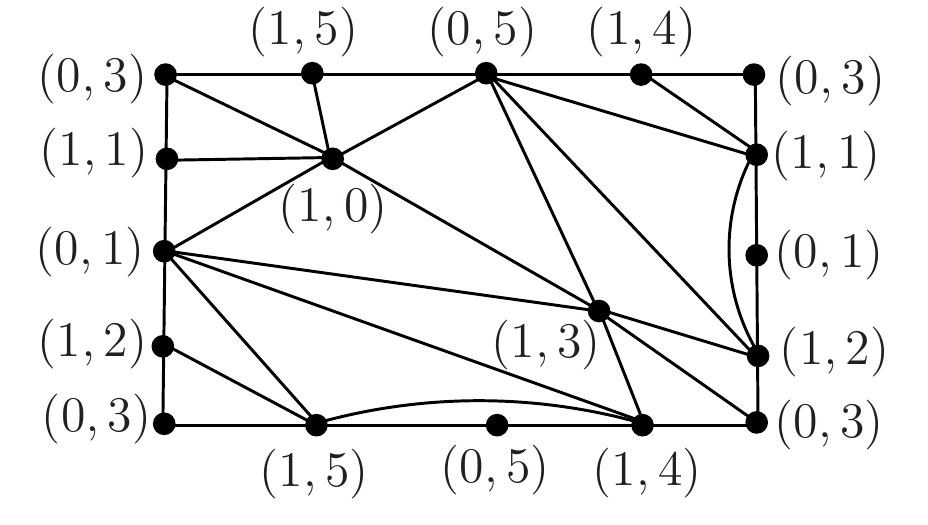}\\
  \caption{An embedding of $\Gamma_{\mathbb{Z}_2 \times \mathbb{Z}_6}$ on $\mathbb{S}_1$.}\label{fn3}
\end{figure}
Now by Lemma \ref{p-group} we see that $\Gamma_G$ has genus one for any group $G$ in (\ref{g1}).

Now we assume that $\gamma(\Gamma_G)=1$. If $G$ is a $p$-group, the desired result
follows from Lemma \ref{p-group}. Thus,
we may assume that $G$ is not a $p$-group. Let $q$ be an odd prime divisor of $|G|$.
If the number of subgroups of order $q$ is not $1$, then by Theorem \ref{Sythm}
$G$ has at least $4$  subgroups of order $q$ and by $(a)$ of Lemma \ref{fact5},
we have a contradiction.
Thus, $G$ has a unique subgroup of order $q$ and thereby,
every Sylow $q$-subgroup of $G$ is cyclic by Lemma \ref{unpg}. Similarly, we can get that
$G$ has a unique Sylow $q$-subgroup. Note that $G$ is not cyclic.
Thus, we may assume that
$G=P\ltimes Q$, where $P$ is a $2$-group and $Q$ is a cyclic group of odd order.
We next prove that $P$ is not cyclic.

Suppose to the contrary that $P$ is cyclic.
Suppose that $|P|= 2$. Then $G$ is dihedral.  By Lemma \ref{planar}, we see that
$Q\ncong \mathbb{Z}_3$ and so $\varphi(|Q|)\ge 4$ and  $G$ has at least $5$ involutions. This implies that $\Gamma_G$ has a subgraph isomorphic to $K_{5,4}$ that
has two partition sets consisting of $5$ involutions and $4$ generators of $Q$,
a contradiction. Suppose now that $P\cong \mathbb{Z}_4$.
Note that $G$ has at least $3$ cyclic subgroups of order $4$. By $(c)$ of Lemma \ref{fact5},
we get $Q\cong \mathbb{Z}_3$. By checking the groups of order $12$,  $G\cong \langle a,b: a^6=b^4=1,b^2=a^3,b^{-1}ab=a^5\rangle$. It is easy to check that  $\gamma(\Gamma_{G})$ has
a subgraph isomorphic to $K_{4,5}$, a contradiction
If $|P|\ge 8$, since $P$ is not normal in $G$, $G$ has at least $3$ Sylow $2$-subgroups and
since $\varphi(|P|)\ge 4$, a contradiction by $(e)$ of Lemma \ref{fact5}.
This means that $P$ is not cyclic.

Note that $V(\Gamma_P)\subseteq V(\Gamma_G)$. Then $\gamma(\Gamma_P)=0$ or $1$ and by Lemmas~ \ref{planar} and \ref{p-group}, $P$ is isomorphic to one of the following groups:
$$
\mathbb{Z}_2^3,~\mathbb{Z}_2\times \mathbb{Z}_4,~D_8,~Q_8,~\mathbb{Z}_2^2.
$$
First by $(d)$ of Lemma \ref{fact5}, we conclude $P\ncong \mathbb{Z}_2^3$.

\medskip
\noindent {\bf Case 1.} $P\cong \mathbb{Z}_2\times \mathbb{Z}_4$.
\medskip

If $P$ is not normal in $G$, then $G$ has at least $4$ cyclic subgroups of order $4$, a contradiction by $(a)$ of Lemma \ref{fact5}.
Thus, $G\cong P\times Q$ and so $G$ has precisely three involutions and at least
two cyclic subgroups of order $4k$ for some odd prime $k$. Considering the generators of
the two cyclic subgroups of order $4k$ and some involution, we have that
$\Gamma_G$ has a subgraph isomorphic to $K_{4,5}$, a contradiction.

\medskip
\noindent {\bf Case 2.} $P\cong D_8$.
\medskip

If $P$ is not normal in $G$, then $G$ has at least $7$ involutions in the union of all
Sylow $2$-subgroups, a contradiction by $(d)$ of Lemma \ref{fact5}.
Therefore, we may assume that $G\cong P\times Q$. Let $g$ be an element of odd order.
Then $G$ has a cyclic subgroup of order $4|g|$, which has at least $4$ generators $\{g_1,\cdots,g_4\}$.
Now it is easy to see that $\Gamma_G$ has a subgraph isomorphic to $K_{4,5}$ that
has two partition sets $4$ involutions and $\{g_1,\cdots,g_4,a\}$ for some involution $a$, a contradiction.

\medskip
\noindent {\bf Case 3.} $P\cong Q_8$.
\medskip

Note that $Q_8$ has $3$ cyclic subgroups of order $4$. By $(a)$ of Lemma \ref{fact5}, we
may assume that $G\cong P\times Q$. So $G$ has at least $3$ cyclic subgroups of order
$4k$ for some odd prime $k$. By  $(e)$ of Lemma \ref{fact5}, a contradiction.

\medskip
\noindent {\bf Case 4.} $P\cong \mathbb{Z}_2^2$.
\medskip

If $P$ is not normal in $G$, then $G$ has at least $7$ involutions, a contradiction by $(d)$ of Lemma \ref{fact5}. Now we assume that  $G\cong \mathbb{Z}_2^2\times Q$. If $Q\cong \mathbb{Z}_3$, then as desired. Thus, we may assume that $|Q|\ge 5$. Then it is easy to see
that $G$ has at least $3$ cyclic subgroups of order $2k$ for some odd number $k\ne 3$. By $(e)$ of Lemma \ref{fact5}, a contradiction.
\qed

\bigskip

{\noindent \em Proof of Theorem \ref{th2}.} Since $\Gamma_{\mathbb{Z}_2\times \mathbb{Z}_4}$ and $\Gamma_{D_8}$ all have some subgraphs isomorphic to $K_{3,3}$, one has that
$\overline{\gamma}(\Gamma_{\mathbb{Z}_2\times \mathbb{Z}_4})\ge 1$ and $\overline{\gamma}(\Gamma_{D_8})\ge 1$. On the other hand, we may embed $\Gamma_{\mathbb{Z}_2\times \mathbb{Z}_4}$ and $\Gamma_{D_8}$ into $\mathbb{N}_1$ as shown in
Figures \ref{fn1} and \ref{fn2}, respectively. So we have $\overline{\gamma}(\Gamma_{\mathbb{Z}_2\times \mathbb{Z}_4})=\overline{\gamma}(\Gamma_{D_8})=1$.
\begin{figure}[hptb]
  \centering
  \includegraphics[width=8cm]{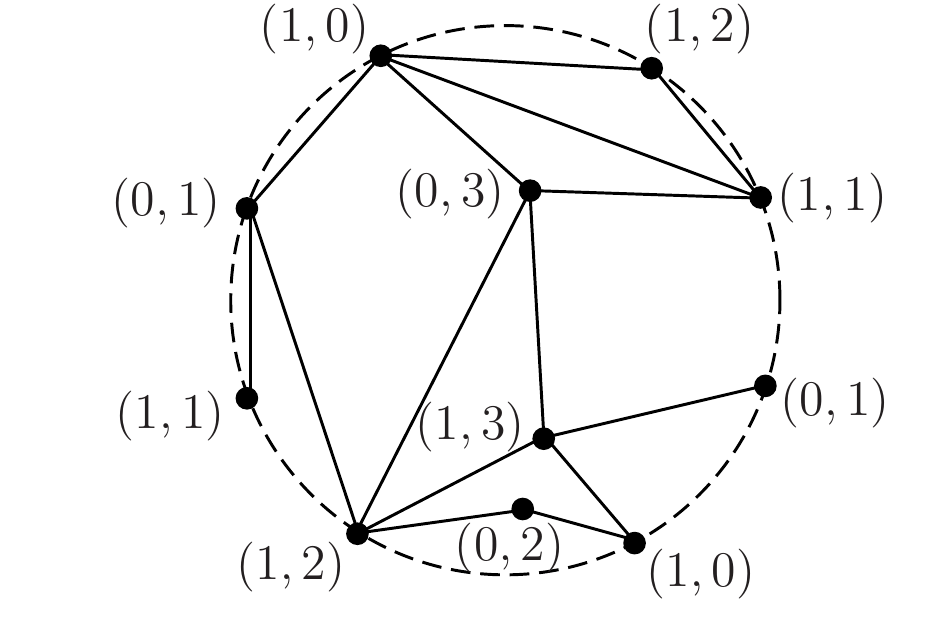}\\
  \caption{An embedding of $\Gamma_{\mathbb{Z}_2\times \mathbb{Z}_4}$ on $\mathbb{N}_1$.}\label{fn1}
\end{figure}
\begin{figure}[hptb]
  \centering
  \includegraphics[width=8cm]{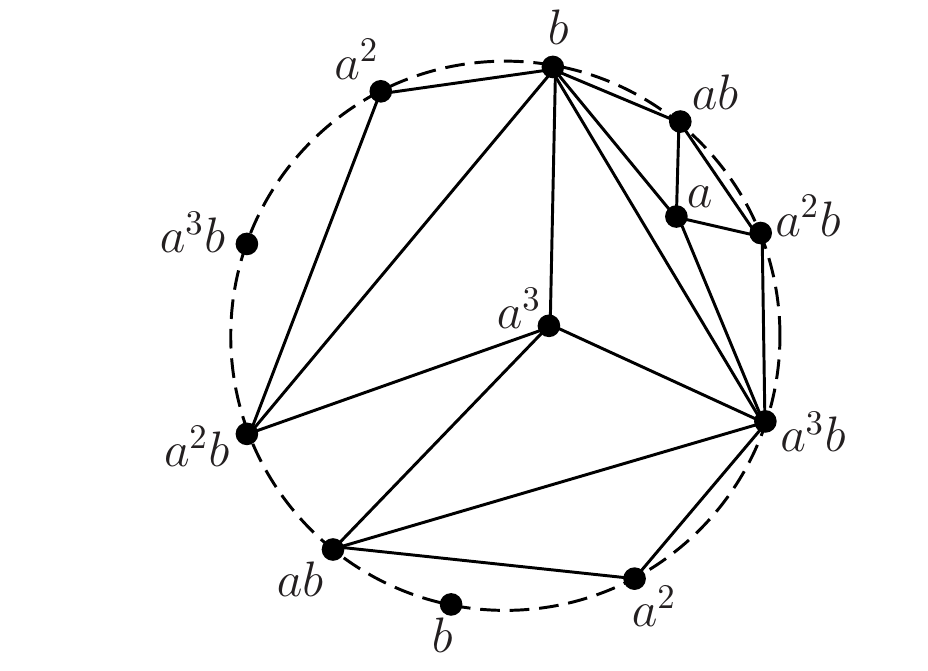}\\
  \caption{An embedding of $\Gamma_{D_8}$ on $\mathbb{N}_1$.}\label{fn2}
\end{figure}

Now we assume that  $\overline{\gamma}(\Gamma_G)=1$.
By Lemma \ref{ccgenus} we see that $\overline{\gamma}(K_{4,4})= 2$
and $\overline{\gamma}(K_{n})\ge 2$ for $n\ge 7$.
Thus, $\Gamma_G$
has no subgraphs isomorphic to $K_{4,4}$ and $K_n$ for $n\ge 7$.
By the proof of Theorem \ref{th1}, it is easy to see that $G$ is one group of (\ref{g1}).
Since $\Gamma_{\mathbb{Z}_3^2}$ has a subgraph isomorphic to $K_{4,4}$ and
$\Gamma_{\mathbb{Z}_2^3}$ has a subgraph isomorphic to $K_{7}$, one has that
$G\cong \mathbb{Z}_2\times \mathbb{Z}_4$, $D_8$ or $\mathbb{Z}_2\times \mathbb{Z}_6$.
In order to complete our proof, we next prove $\overline{\gamma}(\Gamma_{\mathbb{Z}_2\times \mathbb{Z}_6})\ge 2$.

Clearly, $\overline{\gamma}(\Gamma_{\mathbb{Z}_2\times \mathbb{Z}_6})\ge 1$. Suppose for a contradiction that $\overline{\gamma}(\Gamma_{\mathbb{Z}_2\times \mathbb{Z}_6})= 1$.
Note that $|V(\Gamma_{\mathbb{Z}_2\times \mathbb{Z}_6})|=9$ and
$|E(\Gamma_{\mathbb{Z}_2\times \mathbb{Z}_6})|=27$.
Thus, by the Euler characteristic formulas, if $\Gamma_{\mathbb{Z}_2\times \mathbb{Z}_6}$ is embedded into the surface of  nonorientable genus $\overline{\gamma}(\Gamma_{\mathbb{Z}_2\times \mathbb{Z}_6})$, resulting in $f$ faces, then
$$|V(\Gamma_{\mathbb{Z}_2\times \mathbb{Z}_6})|-|E(\Gamma_{\mathbb{Z}_2\times \mathbb{Z}_6})| +f=2-\overline{\gamma}(\Gamma_{\mathbb{Z}_2\times \mathbb{Z}_6}).$$
This implies that
$2|E(\Gamma_{\mathbb{Z}_2\times \mathbb{Z}_6})|\ge 3f$, which is a contradiction as $f=19$.
\qed


\end{document}